\documentclass[12pt]{amsart}

\usepackage{graphics}
\usepackage{amssymb}

\numberwithin{equation}{section}

\newcommand{\fullref}[1]{\ref{#1} on page~\pageref{#1}}

\newcommand{\ndash}{\nobreakdash-\hspace{0pt}}
\newcommand{\Ndash}{\nobreakdash--}

\newcommand{\dd}{{\mathrm{d}}}

\newcommand{\hagrid}{{\calC(M;C_0,C_1)}}

\newtheorem{Thm}{Theorem}[section]

\newtheorem*{Thm*}{Theorem}

\theoremstyle{remark}
\newtheorem{Rem}[Thm]{Remark}
\newtheorem*{Ack}{Acknowledgment}

\theoremstyle{definition}

\newtheorem{Exa}[Thm]{Example}

\newtheorem{Ass}{Assumption}

\newcommand{\braket}[2]{\left\langle{\,{#1}\,,\,{#2}\,}\right\rangle}

\newcommand{\Lie}[2]{{\left[{\,{#1}\,,\,{#2}\,}\right]}}

\newcommand{\Poiss}[2]{\left\{{\,{#1}\,,\,{#2}\,}\right\}}

\newcommand{\bbR}{{\mathbb{R}}}

\newcommand{\de}{\partial}

\newcommand{\calC}{\mathcal{C}}

\newcommand{\calI}{\mathcal{I}}
\newcommand{\calO}{\mathcal{O}}

\newcommand{\frg}{{\mathfrak{g}}}

\newcommand{\frh}{{\mathfrak{h}}}

\def\gpd{\,\lower1pt\hbox{$\longrightarrow$}\hskip-.24in\raise2pt
               \hbox{$\longrightarrow$}\,}

\newcommand\qq{}
\newcommand\cmp[1]{{\qq Commun.\ Math.\ Phys.\ \bf #1}}

\newcommand\anm[1]{{\qq Ann.\ Math.\ \bf #1}}

\newcommand\jdg[1]{{\qq J.\ Diff.\ Geom.\ \bf #1}}

\newcommand\sms[1]{{\qq Selecta Math.\ Soviet.\ \bf #1}}
\newcommand\jmsj[1]{{\qq J. Math.\ Soc.\ Japan \bf #1}}

\begin{document}
\title[Branes in the Poisson sigma model]
{Coisotropic submanifolds in Poisson geometry
and branes in the Poisson sigma model}


\author[A.~S.~Cattaneo]{Alberto~S.~Cattaneo}
\address{Institut f\"ur Mathematik, Universit\"at Z\"urich--Irchel,  
Winterthurerstrasse 190, CH-8057 Z\"urich, Switzerland}  
\email{alberto.cattaneo@math.unizh.ch}

\author[G.~Felder]{Giovanni Felder}
\address{D-MATH, ETH-Zentrum, CH-8092 Z\"urich, Switzerland}
\email{felder@math.ethz.ch}

\thanks{A.~S.~C. acknowledges partial support of SNF Grant No.~20-100029/1}
\thanks{G. F. acknowledges partial support of SNF Grant No.~21-65213.01}

\begin{abstract}
General boundary conditions (``branes'') for the Poisson sigma model are studied.
They turn out to be labeled by coisotropic submanifolds of the given Poisson
manifold. The role played by these boundary conditions both at the classical and
at the perturbative quantum level is discussed. It
turns out to be related at the classical level
to the category of Poisson manifolds with dual pairs as morphisms
and at the perturbative quantum level to the category of associative algebras (deforming
algebras of functions on Poisson manifolds) with bimodules as morphisms.
Possibly singular Poisson manifolds arising from reduction enter naturally into
the picture and, in particular, the construction yields (under certain assumptions)
their deformation quantization.
\end{abstract}

\maketitle

\section{Introduction}\label{intro}
Coisotropic submanifolds play a fundamental role in symplectic geometry as
they describe systems with symmetries (Dirac's ``first-class constraints") and
provide a method to generate new symplectic spaces (``symplectic
reduction").
Their generalizations to Poisson manifolds also carry naturally induced
foliations such that the leaf spaces (``reduced phase spaces'') are again Poisson. They are the
general framework to study symmetries in the Poisson world.
We recall the basic facts about coisotropic submanifolds in Sect.~\ref{coiso}.

The Poisson sigma model \cite{I,SchStr} is a topological field theory defined
in terms of bundle maps from the tangent bundle of a surface to the cotangent
bundle of a given Poisson manifold $M$. Particularly interesting is the case where
the source surface is a disk, which so far has been studied only 
assuming particularly simple boundary conditions (viz., mapping the boundary
to the zero section of $T^*M$); then
the
perturbative path integral expansion  yields \cite{CaFe1} Kontsevich's star product
\cite{Kon} on $M$, while the reduced phase space of the
theory  \cite{CaFeOb} is the symplectic groupoid \cite{K2,W2,Z} of $M$.
A relevant problem concerns other possible boundary conditions and their role.

It turns out that
coisotropic submanifolds of a Poisson manifold label
the possible boundary conditions (``D-branes") of the Poisson sigma model.
Something similar happens in the symplectic context  where
coisotropic submanifolds play an important role as D-branes for the A-model \cite{KO,OP}.

In Sect.~\ref{class} we discuss the
classical Hamiltonian viewpoint. The reduced phase space
of the Poisson sigma model on an interval with boundary conditions labeled by
coisotropic submanifolds $C_0$ and $C_1$ 
is a (possibly singular) symplectic manifold endowed with a Poisson and an anti-Poisson
map to the reduced phase spaces $\underline{C_0}$ and $\underline{C_1}$ of $C_0$ and $C_1$.
This construction yields then 
a ``dual pair"
which is the notion of morphism in a category, whose objects are Poisson manifolds, that
seems to be natural \cite{L} if one has quantization in mind. 

Sect.~\ref{quant}, which can be read independently of Sect.~\ref{class},
deals with the perturbative quantization of the Poisson sigma model
with boundary conditions given by coisotropic submanifolds. We show that locally, under appropriate
assumptions, this construction allows us $i)$  to deformation-quantize
the (possibly singular) Poisson manifolds obtained by reduction from the given
coisotropic submanifolds and $ii)$ to give the space of invariant functions on
the intersection of two coisotropic submanifolds the structure of a bimodule
for the corresponding deformed algebras.
Some examples where the above procedure works are discussed in Sect.~\ref{exa}.

The construction also suggests how to modify Kontsevich's formality map from multivector
fields to multidifferential operators in the presence of a given submanifold, see 
Sect.~\ref{form}. This should be relevant when trying to globalize.

The nonperturbative study (probably beyond our possibilities at the moment) 
looks like a generalization of the Fukaya $A_\infty$\ndash category.

This note is thought of as a short overview of results that will be discussed 
thoroughly elsewhere \cite{CFcoisoclass}. To read Sect.~\ref{class}
the reader is assumed to have had some exposure to \cite{CaFeOb}, while
Sect.~\ref{quant} assumes some familiarity with \cite{CaFe1,Kon}. More advanced remarks,
which have no consequence for the rest of the paper, have been put in footnotes.

\begin{Ack}
We thank Boris Shoikhet, Duco van Straten and Charles Torossian for useful discussions.
We are particularly grateful to James Stasheff for very useful suggestions and for revising
a first version of this paper. 
\end{Ack}

\section{Coisotropic submanifolds}\label{coiso}
A {\sf Poisson manifold} $(M,\pi)$ is a manifold $M$ endowed with a bivector field $\pi$
such that the bracket $\Poiss fg:=\pi(\dd f,\dd g)$ is a Lie bracket on $C^\infty(M)$.
Equivalently, the Poisson bivector field $\pi$ must satisfy $\Lie\pi\pi=0$ where
$\Lie{\ }{\ }$ denotes the 
Schouten--Nijenhuis bracket. In local coordinates, this amounts to the equations
\begin{equation}\label{Jacobi}
\pi^{ij}\,\de_i\pi^{kl}+
\pi^{il}\,\de_i\pi^{jk}+
\pi^{ik}\,\de_i\pi^{lj}
=0.
\end{equation}
The bivector field $\pi$ induces a bundle map $\pi^\sharp\colon T^*M\to TM$ by
\[
\braket{\pi^\sharp(x)\sigma}\tau=\pi(x)(\sigma,\tau), \quad
\forall x\in M,\quad
\forall\sigma,\tau\in T^*_xM,
\]
where $\braket{\ }{\ }$ denotes the canonical pairing.
Some examples of Poisson manifolds are:
\begin{description}
\item[Trivial case] $\pi\equiv0$.
\item[Symplectic case] $(M,\omega)$ is symplectic and $\pi$ is the inverse of $\omega$.
\item[Linear case] $M=\frg^*$, where $\frg$ is a Lie algebra, and the bracket of linear
functions is defined by the Lie bracket. The resulting Poisson structure is usually called
the Kostant--Kirillov Poisson structure.
\end{description}
In general, Poisson manifolds are foliated---by the possibly singular involutive 
distribution\footnote{\label{Moony}The distribution is 
involutive as a consequence of the Jacobi identity
\eqref{Jacobi}. One actually has more structure; viz., $T^*M$ is a Lie algebroid with
$\pi^\sharp$ as its anchor; as for the Lie bracket on its sections, it is enough to define it
on exact $1$\ndash forms for which one sets $\Lie{\dd f}{\dd g}:=\dd\Poiss fg$.
The involutive distribution $\pi^\sharp(T^*M)$ is then the canonical foliation of this Lie 
algebroid.}
$\pi^\sharp(T^*M)$---and each leaf turns out to be symplectic.
In the first example, each point
is a symplectic leaf; in the second example, there is just one symplectic leaf, the whole
manifold; in the third example symplectic leaves are the same as coadjoint orbits
(and have in general varying dimensions).

A submanifold $C$ of a Poisson manifold $(M,\pi)$ is said to be {\sf coisotropic} \cite{W3}
if
$\pi^\sharp(N^*C)\subset TC$, where $N^*C$ denotes the conormal bundle of $C$ (i.e., the subbundle
of $T^*_CM$ consisting of covectors that kill all vectors of $TC$).
It follows from the Jacobi identity for $\pi$ that the 
{\sf characteristic distribution} $\pi^\sharp(N^*C)$ 
on the coisotropic submanifold $C$ is involutive;\footnote{\label{Padfoot}$N^*C$ 
actually turns out
to be a Lie subalgebroid 
of $T^*M$ with Lie algebroid structure as in footnote~\ref{Moony}. 
More precisely, conormal bundles of coisotropic submanifolds
are all possible Lagrangian Lie subalgebroids of $T^*M$ with its canonical symplectic 
structure. If $M$ is integrable, coisotropic submanifolds are also in correspondence with Lagrangian
subgroupoids of the symplectic groupoid of $M$. See \cite{C}.}
the corresponding foliation is called the {\sf characteristic foliation}
and we will denote by $\underline C$ its leaf space which we call 
the {\sf reduced phase space}. 
Its space of ``smooth'' functions
may be defined also when the leaf space is not a smooth manifold
by setting \`a la Whitney $C^\infty(\underline C):=C^\infty(C)^\text{inv}$,
where the superscript denotes the invariant part (a function $f$ on $C$ is invariant
if $X(f)=0$ for all sections $X$ of $\pi^\sharp(N^*C)$).

When $M$ is symplectic, $\pi^\sharp$ yields an isomorphism between $N^*C$ and $T^\perp C$
(the subbundle of $T_CM$ of vectors that are symplectic-orthogonal to all vectors in $TC$).
So we recover the usual definition of coisotropic submanifolds
in the symplectic case: $T^\perp C\subset TC$.

We recall a couple of examples of coisotropic submanifolds. Let $M$ and $N$ be Poisson manifolds
and let $f\colon M\to N$ be a Poisson map
(i.e., a map whose pullback is a morphism of Poisson algebras).
We denote by $\bar N$ the 
Poisson manifold obtained by changing sign to the Poisson structure on $N$.
Then
\begin{enumerate}
\item The graph of $f$ is coisotropic in $M\times\bar N$.
\item The preimage of a symplectic leaf of $N$ is coisotropic in $M$ (when a submanifold).
\end{enumerate}
A particular instance is when $N$ is the dual
of a Lie algebra, in which case
$f$ is an equivariant momentum map. An interesting example, to which we will
return in Sect.~\ref{exa},
is the following:
\begin{Exa}\label{phoenix}
Consider a Lie subalgebra $\frh\stackrel\iota\hookrightarrow\frg$, and set
$M=\frg^*$, $N=\frh^*$ (with Kostant--Kirillov Poisson structure)
and $f=\iota^*$. As $\{0\}$ is a symplectic leaf of $\frh^*$,
we get the coisotropic submanifold $\frh^\perp:=(\iota^*)^{-1}(0)$ (the annihilator of $\frh$)
in $\frg^*$.
\end{Exa}

Let $\calI$ be the ideal  of functions that vanish when
restricted to the submanifold $C$, so $C^\infty(C)=C^\infty(M)/\calI$.
Differentials of elements of $\calI$ yield sections
of $N^*C$. Therefore, we can also characterize coisotropic submanifolds of $M$ as
submanifolds whose vanishing ideal $\calI$ is a Poisson subalgebra (and not just a commutative
subalgebra) of $C^\infty(M)$. 
In Dirac's terminology, a family of functions generating $\calI$ are called first-class
constraints.

Let $N(\calI):=\{g\in C^\infty(M) : \Poiss g\calI\subset\calI\}$ 
be the normalizer of $\calI$. If $\calI$ is a Poisson subalgebra,
so is $N(\calI)$. Moreover, $\calI$ is a Poisson ideal in $N(\calI)$,
so $N(\calI)/\calI$ is a new Poisson algebra. This may easily be recognized as the
algebra $C^\infty(C)^\text{inv}$
of invariant functions on $C$. 
So $\underline C$ is a (possibly singular) Poisson manifold.

Observe that, in the smooth case, the inclusion map $\iota\colon C\to M$ and the
projection $p\colon C\to \underline C$ induce maps of the commutative algebras of functions
that make $C^\infty(C)$ into a bimodule over $C^\infty(\underline C)$ and $C^\infty(M)$.
This clearly works also in the singular case where we have the projection
$\iota^*\colon C^\infty(M)\to C^\infty(M)/\calI$ and the inclusion
$p^*\colon N(\calI)/\calI\to C^\infty(M)/\calI$. 

We may also consider two coisotropic submanifolds $C_0$ and $C_1$. If we denote by
$\underline{C_0\cap C_1}$ the quotient of $C_0\cap C_1$ by the intersection of the characteristic
foliations, we see that $C^\infty(\underline{C_0\cap C_1})$ is a bimodule over the commutative
algebras $C^\infty(\underline{C_0})$ and $C^\infty(\underline{C_1})$. 
(The previous case corresponds to $C_0=C$ and
$C_1=M$.)

The fact that these structures are compatible
with the given Poisson structures gives the bimodule some extra properties that will be better
understood in the following Sections.



\section{Classical Hamiltonian study of the Poisson sigma model}\label{class}
The Poisson sigma model is described at the classical Hamiltonian level by the following data:
$i)$ a weak symplectic structure on an infinite-dimensional manifold (the ``phase space'')
and $ii)$ equations that
select a coisotropic submanifold. As in every topological field theory the 
Hamiltonian is zero and
the characteristic foliation of the coisotropic submanifold has finite codimension (``finitely many
degrees of freedom'').

These data depend on a given Poisson manifold as follows. Let $(M,\pi)$ be a Poisson manifold.
Then the phase space is the cotangent bundle $T^*PM$ of the path space of $M$ (open case) or
the cotangent bundle $T^*LM$ of the loop space of $M$ (closed case) with canonical weak symplectic structure.
These spaces may also be understood as the spaces of bundle maps $TI\to T^*M$ and $TS^1\to T^*M$,
respectively (where $I$ is the interval and $S^1$ the circle).
They may be given a Banach manifold structure by imposing certain conditions (e.g., requiring the base
maps to be differentiable and the fiber maps to be continuous).

An element of these spaces is then a pair $(X,\zeta)$ where $X$ is a (differentiable) map from $I$ or $S^1$ to
$M$ and $\zeta$ is a (continuous) $1$\ndash form taking values in
sections of the pulled-back bundle $X^*T^*M$.
The coisotropic\footnote{To define $\calC(M)$ one just needs a tensor $\pi$. 
One may show however (\cite{SchStr} for the closed and
\cite{CaFeOb} for the open case)
that $\calC(M)$ is coisotropic if{f} $\pi$ is a Poisson bivector field.}
submanifold $\calC(M)$ is defined by the equations
\begin{equation}\label{e:dX}
\dd X+\pi^\sharp(X)\,\zeta=0.
\end{equation}
The characteristic foliation is better described by choosing local coordinates $\{x^I\}_{I=1,\dots,\dim M}$,
so that $X$ and $\zeta$ are locally a set of functions $X^I$ and of $1$\ndash forms $\zeta_I$.
Denoting by $\delta X$ and $\delta\zeta$ the horizontal and vertical components of a vector field
on $T^*PM$, an element of the characteristic distribution 
 is given by
\begin{subequations}\label{e:delta}
\begin{align}
\delta X^I &= \pi^{IJ}(X)\,\beta_J,\label{e:deltaX}\\
\delta\zeta_I &= 
-\dd\beta_I -\de_I\pi^{JK}\,\zeta_J\,\beta_K,\label{e:deltazeta}
\end{align}
\end{subequations}
where $\beta$ is a (differentiable) section of $X^*T^*M$ that, in the open case,
is required to vanish on the boundary.
The reduced phase space\footnote{\label{f:morphLA} Using the language of Lie algebroids,
one may also give the following interpretation \cite{S,CF}: 
Elements of $\calC(M)$ are precisely those bundle maps
that are also morphisms of Lie algebroids, where the tangent bundles are given the canonical Lie algebroid
structure and $T^*M$ the one induced from the Poisson structure (see footnote~\ref{Moony}). 
Elements of $\underline{\calC(M)}$
are then morphisms of Lie algebroids modulo ``homotopy.'' In the closed case, we say that
two morphisms $\gamma_0,\gamma_1\colon TS^1\to T^*M$ are homotopic if
there exists a morphism of Lie algebroids $T(S^1\times [0,1])\to T^*M$ such that
its restriction to $T(S^1\times\{u\})$ is $\gamma_u$, $u=0,1$.
In the open case, beside the obvious replacement of $S^1$ by $I$, we put the extra condition
that the restriction of the morphism to $T(\de I\times[0,1])$ is the zero bundle map
(or, in other words, a morphism to the rank-zero Lie algebroid over $M$ regarded as a Lie subalgebroid
of $T^*M$.)}
$\underline{\calC(M)}$ has particularly interesting properties
in the open case (where it is shown to be the possibly singular, source-simply-connected
symplectic groupoid of $M$ \cite{CaFeOb}).

{}From now we will consider only the open case and look for possible boundary conditions. Given
two submanifolds $C_0$ and $C_1$ of $M$, we define $\hagrid$ to be the submanifold of $\calC(M)$ 
where the base maps are paths connecting $C_0$ to $C_1$ (with this new notation we have, in particular,
$\calC(M)=\calC(M;M,M)$). We have then \cite{CFcoisoclass} the following:
\begin{Thm}
Assume that all pairs of points of the two coisotropic submanifolds can be connected by base paths of solutions of \eqref{e:dX}. Then 
$\hagrid$ is coisotropic in
$T^*PM$ if{f} $C_0$ and $C_1$ are coisotropic in $M$. 
\end{Thm}
The characteristic distribution is again given 
by \eqref{e:delta} but with the condition that $\beta$ on the boundary $\de I=\{0,1\}$
is an element of $N^*_{X(u)}C_u$, $u=0,1$. (The previous case is obtained by setting $C_0=C_1=M$
and observing that $N^*M$ has rank zero.) Observe that the coisotropy condition on $C_t$, $t=0,1$,
ensures that $\delta X(t)$ is tangent to $C_t$, as required by the boundary conditions.

The characteristic foliation\footnote{The leaf space $\underline\hagrid$ may also be
defined as in footnote~\ref{f:morphLA} as the quotient of the space of Lie algebroid
morphisms $TI\to T^*M$ by homotopies. The morphisms are however now required to have base maps
connecting $C_0$ to $C_1$, and homotopies must satisfy the condition that
the restriction to $T(t\times[0,1])$, $t=0,1$, is a morphism of Lie algebroids with range
$N^*C_t$.}
on $\hagrid$
may move the
endpoints of the base maps but only along the characteristic foliations of $C_0$ and $C_1$.
Thus, the evaluation maps at $0$ and $1$ descend to the quotients and define maps 
$J_u\colon\underline\hagrid\to\underline{C_u}$, $u=0,1$.
Observe that, when smooth, $\underline\hagrid$ is endowed with a symplectic structure while
$\underline{C_0}$ and $\underline{C_1}$ are endowed with Poisson structures. It is then possible to prove
\cite{CFcoisoclass} the following:
\begin{Thm}
$J_0$ and $J_1$ are a Poisson and an anti-Poisson map, respectively, and the $J_0$\ndash fibers
are symplectically orthogonal to the $J_1$\ndash fibers (so pullbacks of functions
via $J_0$ and $J_1$ Poisson commute).
\end{Thm}
Thus, using the terminology of \cite{K,W} (see also \cite{BW,L} and references therein),
$\underline{C_0}\stackrel{J_0}\gets\underline\hagrid\stackrel{J_1}\to\underline{C_1}$
is a (possibly singular) {\sf dual pair}.
Observe \cite{L1}
that dual pairs are the morphisms of a category in which Poisson manifolds are the objects
(the composition of the dual pairs $S$ and $S'$ which have the same Poisson manifold $P$
as target and source, respectively, is obtained by symplectic reduction observing that
$S\times_PS'$ is coisotropic in $S\times S'$). 
This structure suggests, given a Poisson manifold, 
to define a category\footnote{There is also another category, actually a groupoid, associated
to the Poisson sigma model with boundary: its objects are points in the reduced coisotropic submanifolds and
the morphisms between $[x_0]\in\underline{C_0}$ and $[x_1]\in\underline{C_1}$ are the elements of
$\underline\hagrid$ 
with $J_t([(X,\zeta)])=[x_t]$. Composition is obtained by gluing, and the inverse by reversing
$I$. The symplectic groupoid of $M$ is then a subgroupoid of this.}
whose objects are the leaf spaces of its coisotropic submanifolds and whose
morphisms are generated by the dual pairs obtained above.

In Sect.~\ref{quant} we will see (cf.\ Thm.~\ref{Luna}) that the corresponding
quantum category (in the context of deformation quantization) consists of associative algebras
with bimodules as morphisms.

\section{Perturbative quantization of the Poisson sigma model}\label{quant}
\subsection{Classical action functional and symmetries}
In the path integral quantization of the Poisson sigma model,
one starts from a classical action functional $S$,
a function on the space of bundle maps $T\Sigma\to T^*M$
from the tangent bundle of a surface $\Sigma$ to the
cotangent bundle of the Poisson manifold $M$. Such a bundle
map $\hat X$ consists of a base map $X\colon\Sigma\to M$ and a linear
map $\eta$ for each fiber, which may be thought of 
as a 1-form $\eta\in\Omega^1(\Sigma,X^*T^*M)$ on $\Sigma$
with values in the pull-back of the cotangent bundle.
The action functional is then \cite{I, SchStr}
 $S(X,\eta)=\int_\Sigma(\langle\eta,dX\rangle+
\frac12\langle\pi\circ X,\eta\land\eta\rangle)$.
In the case of interest to us where $\Sigma$ has a boundary,
it is natural to consider the action functional 
with boundary conditions imposing that $\hat X$ maps the
tangent bundle $T\partial\Sigma$ of the boundary to
the conormal bundle $N^*C$ of a submanifold $C$.
 With these boundary conditions, 
the Euler--Lagrange equations are 
differential equations without any
boundary term, since the boundary term coming from
integration by parts is 
\begin{equation}\label{e-Septimus}
\int_{\partial\Sigma}\langle\eta,\delta X\rangle,
\end{equation}
which vanishes for any variation $\delta X$ of the
base map. 

If $C$ is the whole of $M$, this boundary
condition is the one considered in \cite{CaFe1} and leads,
in case $\Sigma$ is a disk,
to the construction of the Kontsevich formula for
deformation quantization of $M$. In this case
there are no conditions on the base map $X$ and
$\eta$ is assumed to vanish on vectors tangent to the
boundary of $\Sigma$. 

If $C$ is a coisotropic submanifold, 
the boundary conditions for gauge transformations 
of \cite{CaFe1} can be generalized to this case. An
infinitesimal gauge transformation at $\hat X$ is 
parametrized by a section 
$c\in\Gamma(\Sigma,X^*T^*M)$ restricting
to the boundary to a section of $X^*N^*C$. 
The action functional is invariant under such a
gauge transformation if $C$ is coisotropic.
Indeed, the calculation of the variation of the action
of \cite{I,SchStr}, done for closed $\Sigma$, shows
that $S$ is invariant up to the boundary term
\eqref{e-Septimus}. The infinitesimal
gauge
variation of the base map is $\delta X=\pi^\sharp{}c{}$,
so that the boundary term vanishes if $C$ is coisotropic.

{}From now on we restrict our attention to the case where
$\Sigma$ is a disk.

\subsection{Batalin--Vilkovisky quantization}
The quantization of the Poisson sigma model with
boundary conditions is given by path integrals
$\int \exp(i S/\hbar) \mathcal O \dd\hat X$ over the
space of bundle maps $\hat X=(X,\eta)$ 
obeying the boundary conditions. 
The observables $\mathcal O$ are gauge invariant
functions on this space. A class of observables 
of particular interest is given by evaluating functions
on $M$ at the image by $X$ of the points of the boundary:
$\mathcal O=\prod_{i=1}^k f_i(X(p_i))$, 
$p_i\in\partial \Sigma$. The condition of gauge invariance
is then $\pi( \dd f_i,{}c{})=0$ for ${}c{}\in N^*C$, i.e.,
$f_i\in N(\mathcal I)$. Since only the value of
$f_i$ on $C$ matters we may take $f_i\in N(\mathcal I)/\mathcal I=C^\infty(\underline C)$. 
Equivalently, the functions $f_i$ are functions
on $C$ which are constant on the leaves of the 
foliation.

This reasoning and the results of \cite{CaFe1}, where
the case $C=M$ was considered, suggest that the
Batalin--Vilkovisky
 perturbative calculation of the path integral should
yield an associative product on $C^\infty(\underline C)$
obtained by picking three distinct points $p,q,r$
on the boundary of the disk $\Sigma$ and setting
\begin{equation}\label{e-Lucrezia}
(f\star g)(x)=\int_{X(r)=x}e^{\frac i\hbar 
S(X,\eta)}f(X(p))g(X(q))\dd X \dd\eta,
\end{equation}
$f,g\in C^\infty(\underline C)$.
The Batalin--Vilkovisky procedure gives a way to make
sense (as a formal power series in $\hbar$)
of this integral by deforming the integration
domain to a Lagrangian submanifold of the odd
symplectic $Q$\ndash manifold of maps $\Pi T\Sigma\to
\Pi T^*M$, see \cite{CaFe1,CaFeAKSZ}, giving a 
version of an AKSZ model \cite{AKSZ}. 
This essentially amounts to replacing $(X,\eta)$
by superfields $(\mathbf X,\boldsymbol\eta)$, where
$\mathbf X$ is a map $\Pi T\Sigma\to M$ to the base
and $\boldsymbol\eta$ is a section of the pull-back
$\mathbf X^*\Pi T^*M$. 
The action functional is $S=S_0+S_\pi$, where $S_0
=\int_{\Pi T\Sigma} \langle\boldsymbol\eta,D \mathbf X\rangle
\mu$
and for any multivector field $\alpha$, $S_\alpha
=\int_{\Pi T\Sigma}\langle\alpha\circ \mathbf X,\boldsymbol\eta\land
\cdots\land\boldsymbol\eta\rangle\mu$. Here $\mu$ is the canonical
volume form on $\Pi T\Sigma$ and $D$ is induced by the
de~Rham differential on $C^\infty(\Pi T\Sigma)=
\Omega^\cdot(\Sigma)$.

The boundary conditions for the case $C=M$ were discussed
in \cite{CaFe1,CaFeAKSZ}. Similar arguments apply here.
The result is that the classical master equation 
$\{S,S\}=0$ is obeyed if 
the boundary conditions are that $(\mathbf X,\boldsymbol\eta)$ 
restricts on the boundary to a map $\Pi T\partial\Sigma
\to\Pi N^*C$  for a coisotropic submanifold $C$.\footnote{In the AKSZ formulation,
the possible boundary conditions are discussed in \cite{CaFeAKSZ}. If 
the source supermanifold is of the form $\Pi T\Sigma$ ($\Sigma$ a manifold with boundary)
and  the target supermanifold $Y$
has a $QP$\ndash structure defined by an odd symplectic form $\omega=\dd\theta$ and
a solution $s$ of the master equation, then the boundary conditions are labeled by Lagrangian
submanifolds of $Y$ where both $\theta$ and $s$ restrict to zero (and, given such a Lagrangian
submanifold $L$, one requires maps $\Pi T\Sigma\to Y$
to restrict on the boundary to maps $\Pi T\de\Sigma\to L$). 
In the present case $Y=\Pi T^*M$, $\theta$ is the canonical $1$\ndash form
$\braket p{\dd x}$ and, given a Poisson bivector field $\pi$, we set $s=\braket p{\pi^\sharp(x)p}/2$
(we denote by $x$ coordinates on $M$ and by $p$ coordinates on the fiber).
So Lagrangian submanifolds with the above properties are the same as odd conormal bundles of coisotropic
submanifolds of $M$. 

Similar boundary
conditions for the A-model are proposed in \cite{OP} where $M$
is symplectic and $N^*C$ is replaced by $T^\perp C$.} 
Here the curly bracket (the BV bracket) 
denotes the Poisson bracket associated
to the odd symplectic structure. Indeed, we have
in general $\{S_\alpha,S_\beta\}=S_{[\alpha,\beta]}$,
 so that 
$\{S_\pi,S_\pi\}=0$ for Poisson bivector fields $\pi$. The bracket 
with $S_0$ involve a boundary term from integration
by parts. With our boundary conditions,
$\{S_0,S_0\}$ vanishes (for any $C$) and
$\{S_0,S_\pi\}$ 
vanishes for $C$ coisotropic as the boundary
term is proportional to $\int_{\partial\Sigma}
\langle \pi\circ\mathbf X,\boldsymbol\eta\land\boldsymbol\eta\rangle$.
The observables $\mathcal O$ are then cocycles for the BV differential
$\{S,\ \}$.\footnote{\label{Prongs}As 
observed in footnote~\ref{Padfoot},
the conormal bundle $N^*C$ of $C$ is a Lagrangian Lie subalgebroid of $T^*M$, so $\Pi N^*C$
is a Lagrangian submanifold of $\Pi T^*M$.
One may then define more general boundary observables associated to elements of the
corresponding Lie algebroid cohomology (invariant functions being the case of degree zero).

In fact, let $V$ be a representative of a Lie algebroid cohomology class of degree $k$. In particular,
$V$ is a section of the
$k$th exterior power of the normal bundle $NC=T_CM/TC$.  With our choice of
coordinates, we may write $V=V^{\mu_1\dots\mu_k}\de_{\mu_1}\land\dots\land\de_{\mu_k}$.
To it, we associate the functional
\[
\mathbf V:= V(\mathbf X)^{\mu_1\dots\mu_k}\, \boldsymbol\eta_{\mu_1}\land\dots\land\boldsymbol\eta_{\mu_k}.
\]
We then get observables either by evaluating $\mathbf V$ at a point $p\in\de\Sigma$,
\[
\calO^0_{V} := \mathbf V(p) = V(X(p))^{\mu_1\dots\mu_k}\,c_{\mu_1}(p)\cdots c_{\mu_k}(p),
\]
or by integrating it on the whole boundary, $\calO^1_{V} := \int_{\de\Sigma} \mathbf V$.
It turns out that $\calO^0_V$ and $\calO^1_V$ are BV closed observables (of degree $k$ and $k-1$ respectively)
and that their BV cohomology classes are independent of the choices above.
}

More generally one may take $k$ coisotropic submanifolds
$C_0,\dots,C_{k-1}$ and consider $\Sigma$ to be a
disk whose boundary is partitioned into  $k$  intervals
with the boundary condition that
$\hat X$ maps the tangent bundle of the $i$th interval $I_i$  
to the conormal bundle $N^*C_i$. The gauge parameter
${}c{}$ maps $I_i$ to $N^*C_i$. Gauge invariant 
observables
are obtained by evaluating functions in $C^\infty(\underline {C_i})$ at 
the image of points in the interior of $I_i$ or
functions in $C^\infty(\underline{C_{i}\cap C_{i+1}})$
evaluated at the point separating two neighboring
intervals $I_i$ and $I_{i+1}$, $i=0,\dots, k-2$. The point $r$ separating
$I_{k_1}$ and $I_0$ is used to select a classical solution by the condition
$X(r)=x$. 

\subsection{Deformation of bimodule structures}
In the next to simplest case $k=2$, we then have two
submanifolds $C_0$, $C_1$ and divide the boundary of
the disk $\Sigma$ into two intervals $I_0,I_1$ whose
common boundary points are two points $p,q\in \partial \Sigma$. Considering path integrals with these boundary conditions and the condition that $X(q)=x\in C_0\cap C_1$ we
obtain various products between functions in
$C^\infty(\underline {C_i})$ and $C^\infty(\underline
{C_0\cap C_1})$, depending on the points on $\partial\Sigma$ 
at which we evaluate the functions.
The associativity of these products are then expected to
give a deformation of the $C^\infty(\underline {C_0})$--%
$C^\infty(\underline {C_1})$\ndash bimodule structure of
$C^\infty(\underline{C_0\cap C_1})$, where the deformation of the
product on $C^\infty(C_i)$ is obtained from the case $k=1$ 
considered above. 
Observe that associative algebras with bimodules as morphisms form a
category which is in some sense the quantization of the category of dual
pairs described at the end of Sect.~\ref{class}.

Of course these semiclassical
statements are expected to receive
quantum corrections and should not be expected to
hold without some additional assumptions. In fact
we consider here only very simple situations in which
the perturbative expansion can be computed and the
statements can be checked rigorously at the level of finite-%
dimensional Feynman integrals. 

\subsection{Feynman expansion}
We start from the case of one coisotropic submanifold and 
consider the case where $M$ is an open subset of
$\mathbb R^n$
with coordinates $x^1,\dots,x^n$ and the
submanifold $C$ is given
by the equations 
\begin{equation}\label{e-PeterWalsh}
x^\mu=0, \qquad \mu=m+1,\dots, n,
\end{equation}
The tangent space to a point on $C$ is then spanned
by $\partial/\partial x^i$, $i=1,\dots,m$, and the
conormal bundle by $dx^\mu$, $\mu=m+1,\dots, n$.
With the convention that lower case Latin indices
run over $\{1,\dots,m\}$ and Greek indices over
$\{m+1,\dots,n\}$, the condition of coisotropy is
then the condition that 
\[
\pi^{\mu\nu}(x^1,\dots,x^m,0,\dots,0)=0,
\]
for the components of the tensor $\pi$. The coordinate
functions $x^\mu$ are a system of generators for the
ideal of $C$ and the characteristic foliation is spanned
by the vector fields 
$E^\mu=\sum_{i=1}^n\pi^{\mu i}\partial_i$ on $C$.
The invariant functions on $C$ are solutions of
\begin{equation}\label{e-Clarissa}
E^\mu(f)=\sum_{i=1}^m\pi^{\mu i}\partial_if=0.
\end{equation}
This condition will be modified by terms of higher order
in $\epsilon$.

The boundary 
conditions for the superfield are then $\mathbf X^\mu=0$,
$\boldsymbol\eta_i=0$ on $\Pi T\Sigma$.

The evaluation of the integral 
\eqref{e-Lucrezia} in a power series in $\epsilon$
along the lines of \cite{CaFe1} leads to a
modification of the Kontsevich formulas of \cite{Kon}.
They can be written as follows
\begin{equation}\label{e-SallySeton}
f\star g=fg+\sum_{k=1}^\infty\frac{\epsilon^k}{k!}
\sum_{\Gamma\in G_{k,2}} w_\Gamma B_\Gamma(f,g),
\qquad f,g\in C^\infty(C).
\end{equation}
The sum is over admissible graphs $\Gamma$ of order $k$,
to which are associated a 
{\em weight} $w_\Gamma\in\mathbb R$
and a {\em bidifferential operator} $B_\Gamma$. The deformation
parameter is $\epsilon=i\hbar$.

\begin{figure}[b]
\begin{picture}(0,200) (100,20)
{\includegraphics{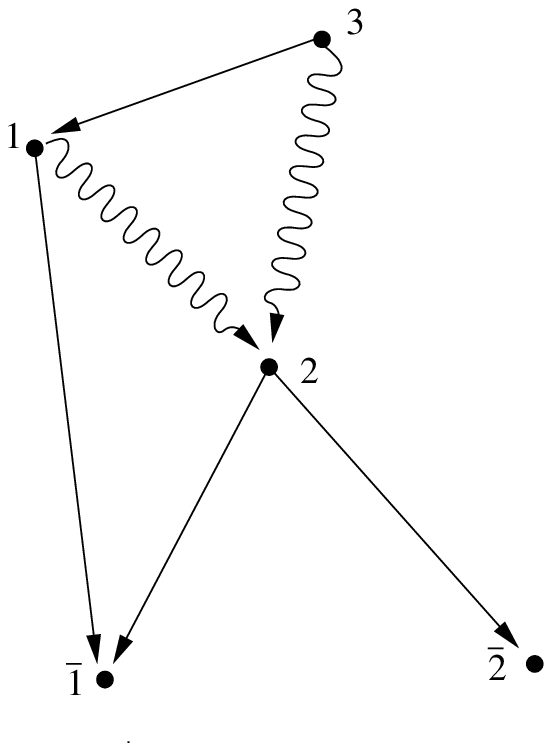}}
\end{picture}
\caption{A simple admissible graph}\label{fig1}
\end{figure}

An admissible graph in $G_{k,2}$ has $k$ vertices  
$1,\dots,k$ of the first type, and 2 vertices $\bar 1,
\bar 2$ of the second type. 
The edges are oriented and come in two types,
say straight and wavy. There are exactly two edges emerging
from each of the vertices of the first type
and none from vertices of the second type.
An ordering of the edges emerging from each vertex is
given.
Each edge may land at any vertex except at 
the one it emerges
from. A simple example of such a graph is given
in Fig.~\ref{fig1}.

The bidifferential operator associated to $\Gamma$
is obtained by the following rule: to each vertex 
of the first type we associate a component of $\pi$ and
to the vertices of the second type
 we associate the functions
$f$ and $g$. To an edge from a vertex
$a$ to a vertex $b$ we associate a partial derivative
acting on the object associated to $b$ with respect
to the variable with the same index as the corresponding
index  of the component of $\pi$ associated to $a$.
Then we take the product and sum over Latin indices for
straight lines and over Greek indices for wavy lines.
Finally we evaluate the result at a point $x\in C$.
For example $\Gamma$ in Fig.~\ref{fig1} gives the
bidifferential operator
\[
\partial_l                   \pi^{i\lambda}
\partial_\lambda\partial_\mu \pi^{jk}
                             \pi^{l\mu}
\partial_i\partial_j f
\partial_k g.
\]
The sum over $\{1,\dots,m\}$ for
repeated Latin indices
and over $\{m+1,\dots,n\}$ for repeated
Greek indices is
understood.

The weight of $\Gamma$ is 
\[
w_\Gamma=\frac1{(2\pi)^{2k}}\int_{C^+_{k,2}}
\;\prod_{\mathrm{edges\; e}}\dd\phi_e.
\]
The integral is over the configuration space $C^+_{k,2}$
of $k$ distinct points $z_i$ in the upper half plane and
two ordered points $z_{\bar 1}<z_{\bar 2}$
on the real axis, modulo dilations
and translation in the real direction. The differential
form $\dd\phi_e$ associated to an edge $e$ going from
$a$ to $b$ is $\dd\phi(z_a,z_b)$ if the edge is straight
and is $\dd\phi(z_b,z_a)$ if it is wavy. Here $\dd\phi(z,w)$
is the differential of the Kontsevich angle function
\[
\phi(z,w)=\frac1{2i}\log\,
\frac{(z-w)(z-\bar w)}{(\bar z-w)(\bar z-\bar w)}
=\mathrm{arg}(z-w)+\mathrm{arg}(z-\bar w).
\]
The ordering of factors in the product of $\dd\phi_e$ is
obtained from
the ordering of vertices and the given
ordering of edges emerging from each vertex.

The fact that wavy lines correspond to $\dd\phi(z_b,z_a)$ 
rather than to
$\dd\phi(z_a,z_b)$ and the fact that the result of the action of
the bidifferential operators are evaluated 
at a point $x\in C$ are
the only places where the formula differs from 
Kontsevich's (the case
$C=M$). For the readers familiar with \cite{CaFe1} 
we add that  the (super-)propagators 
$\langle\boldsymbol\eta_I(z)\boldsymbol\xi^J(w)\rangle=
\delta_I^J P_I(z,w)$ in the Feynman 
perturbative expansion around the constant classical
solution $\mathbf X(z)=x$, $\boldsymbol\eta=0$ 
($\boldsymbol\xi^I=\boldsymbol X^I-x^I$) 
differ for $I=\mu\in \{m+1,\dots,n\}$
from the ones in the case $C=M$ by the
boundary condition, which is that it vanishes when $w$
rather than $z$ is restricted to the boundary. So we have
$P_i(z,w)=\dd\phi(z,w)$ as in the case $C=M$ but $P_\mu(z,w)=
\dd\phi(w,z)$. 

Note that as the differential associated to a wavy edge
vanishes if it points to $\bar1$ or $\bar2$, the functions
$f$, $g$ are differentiated only in the tangential 
directions $\partial_j$. Therefore the bidifferential operators
$B_\Gamma$ are well-defined on functions on $C$.

\subsection{Stokes' theorem and associativity}\label{s-43}
As in \cite{Kon}, the main tool to prove properties of
the product is Stokes' theorem on a 
compactification $\bar C^+_{k,m}$
of configuration spaces $C^+_{k,m}$ of $k$ distinct
points in the upper half-plane and $m$ ordered points
on the real axis modulo translations and dilations. 
For example, the
associativity of the Kontsevich product \eqref{e-SallySeton} (the case $C=M$)
is proven by evaluating the integral of the differential
of a closed form (which of course vanishes) on 
$\bar C^+_{k,3}$ with Stokes' theorem. The sum of contributions 
of the faces (pieces of the boundary) yield associativity
identities.

 The same calculation can be applied to the case
of general $C$ of the type \eqref{e-PeterWalsh}, but there
is an important difference: the contribution
from some faces (the faces with a ``bad edge'') does
does not vanish a priori. These are faces containing
limiting configurations where a subset of the points
in the upper half-plane approach a  
point on the real axis.  
These faces produce corrections to the associativity
involving more general objects expressed in terms of
graphs, which we proceed to describe.

For vector fields $\xi,\eta$ on $M$ 
introduce a differential
operator $A(\xi)$ on $C^\infty(C)[[\epsilon]]$ by
\[
A(\xi)f=\xi f+\sum_{k=1}^\infty\frac{\epsilon^k}{k!}
\sum_{\Gamma\in G_{k+1,1}} w_\Gamma B_\Gamma(\xi)f
\]
and a function
\[
F(\xi,\eta)=\sum_{k=0}^\infty\frac{\epsilon^k}{k!}
\sum_{\Gamma\in G_{k+2,0}}w_\Gamma B_\Gamma(\xi,\eta)
\in C^\infty(C)[[\epsilon]].
\]
The definitions of $w_\Gamma$, $B_\Gamma$
are the same as for $G_{k,2}$ except
that graphs in 
$G_{k+1,1}$ have one additional vertex of the first type
associated with $\xi$
from which one line emerges and just one vertex of the
second type;  graphs in
$G_{k+2,0}$ have two additional vertices of the first type
associated to $\xi$ and $\eta$ and none of the second 
type. In the case $C=M$ these objects were introduced in \cite{CaFeTo}
to construct global star-products on manifolds.

{}From now on, we make the following

\medskip

\begin{Ass}\label{a-Elisabeth}
 $F(E^\mu,E^\nu)=0$ for $m+1\leq\mu,\nu\leq n$.
\end{Ass}

\medskip

This assumption is verified in a number of examples, as
we discuss below. It appears that it is possible to remove
this assumption at the cost of introducing a recursive
correction procedure. This will be discussed elsewhere \cite{CFcoisoclass}.

The quantum version of the algebra of invariant
functions on $C$ is defined to be the $\mathbb R[[\epsilon]]$\ndash module
\[
C_\epsilon^\infty(\underline C)=\{
f\in C^\infty(C)[[\epsilon]]:
A(E^\mu)f=0\}.
\]

\begin{Thm}\label{t-LadyBruton}
Under Assumption~\ref{a-Elisabeth}
the product \eqref{e-SallySeton}
restricts to an
associative product on $C_\epsilon^\infty
(\underline C)$
\end{Thm}

The proof is similar to Kontsevich's proof of associativity
of his star-product and is based on Stokes' theorem.
In this case new boundary components give potentially non-trivial
contributions due to the fact that the 1-form associated
to wavy edges does not vanish as the first argument 
approaches the real axis. These contributions vanish
under Assumption~\ref{a-Elisabeth} and the condition
defining $C_\epsilon^\infty(\underline C)$. Details
will appear elsewhere \cite{CFcoisoclass}.

\begin{Rem}
In general, $(C_\epsilon^\infty(\underline C),\star)$
is not a deformation of $C^\infty(\underline C)$. What
we have is that the map
\[
p\colon f_0+\epsilon f_1+\epsilon^2f_2+\cdots\mapsto f_0
\]
is a ring homomorphism 
$(C_\epsilon^\infty(\underline C),\star)
\to (C^\infty(\underline C),\cdot)$ with the property
that $p((f\star g-g\star f)/\epsilon)=\{f,g\}$.
It would be interesting to characterize the image of
this homomorphism.
\end{Rem}

\subsection{The case of two coisotropic submanifolds:
bimodules}\label{s-44}
The above calculation may be extended to the case of
an arbitrary number of coisotropic submanifolds. We
discuss here the simplest case of two cleanly intersecting
submanifolds $C_0,C_1$ (one says that the intersection of $C_0$ and $C_1$
is clean if $C_0\cap C_1$ is also a submanifold and 
$T(C_0\cap C_1)=TC_0\cap TC_1$).
Again we consider path integrals of the type \eqref{e-Lucrezia} and evaluate the product at a point 
$x\in C_0\cap C_1$. The circle is partitioned into two parts which are
sent to the two coisotropic submanifolds. It is convenient
to map the disk to the first quadrant 
$\mathrm{Re}\,z\geq 0, \mathrm{Im}\,z\geq 0$. 
The parts of the boundary sent
to $C_0$, $C_1$ are the positive imaginary and real
axes, respectively. The point $r$ which is sent to
the point at which we evaluate the product in
\eqref{e-Lucrezia} is the point
at infinity. We then have the option of putting the
remaining points at $0$ or on the real or imaginary axis,
obtaining various products. 

Specifically, let us consider the case where $M$
is an open subset of $\mathbb R^n$
 and suppose that $C_q$, $q=0,1$, is given
by the equations
\[
x^\mu=0,\qquad \mu\in I_q^c,
\]
for subsets $I_0, I_1$ of $\{1,\dots,n\}$, with
complements $I_0^c,I_1^c$. Then
$x^i$, $i\in I_q$ form a coordinate system for $C_q$
and the intersection $C_0\cap C_1$ has coordinates
$x^i$, $i\in I_0\cap I_1$. Such a choice of coordinates
is possible in the neighborhood of a point of
clean intersection.

We suppose that Assumption~\ref{a-Elisabeth} holds
for both $C_0$ and $C_1$ and have therefore two
algebras $C^\infty_\epsilon(\underline {C_0})$, 
$C^\infty_\epsilon(\underline {C_1})$. The 
evaluation of the path integral gives a
$C^\infty_\epsilon(\underline {C_0})$--%
$C^\infty_\epsilon(\underline {C_1})$\ndash bimodule
$C^\infty_\epsilon(\underline{C_0\cap C_1})$.
The construction is in terms of sums
over graphs and goes as follows.

The set of admissible graphs $G_{k,2}$ consists
in this case of graphs with $k$ vertices 
$1,\dots,k$ of the first kind and $2$ vertices
$\bar 1,\bar 2$ of the second kind. The
rules are as before except that there are four
types of vertices $++$, $+-$, $-+$, $--$,
rather than just two. To each such graph $\Gamma$ one
associates a bidifferential operator 
$B_\Gamma(f,g)$. The rules are
the same as in the case of one submanifold, 
the only difference being
the range of summation of the indices associated
to the edges: an edge of type $++$, $+-$, $-+$, $--$
indicates a summation over 
$I_0\cap I_1$,
$I_0\cap I_1^c$,
$I_0^c\cap I_1$,
$I_0^c\cap I_1^c$, respectively.
We also consider graphs in $G_{k+1,1}$ with one additional
vertex with one outgoing edge and one vertex of the
second type. They give rise to differential operators
$B_\Gamma(\xi)$ depending of a vector field $\xi$.

The weight $w_\Gamma$ of a graph $\Gamma$
is obtained by integrating 
the product of one-forms associated to edges over 
configuration spaces.
The one-forms $\dd\phi_{\sigma\tau}(z,w)$, $\sigma,\tau=\pm1$,
corresponding to the different kinds
of edges are obtained from the Euclidean angle
function $\phi_e(z,w)=\mathrm{arg}(z-w)$ by reflection:
\[
\phi_{\sigma\tau}(z,w)=
            \phi_e(z,      w)
+\sigma     \phi_e(z, \bar w)
+\tau       \phi_e(z,-\bar w)
+\sigma\tau \phi_e(z,-     w).
\]
If $\Gamma\in G_{k+1,1}$ the integration is over
the configuration space of $k+1$ points in the
first quadrant modulo dilations. The 
differential operator $A(\xi)=\sum_{\Gamma\in G_{k+1,1}}
w_\Gamma B_\Gamma(\xi)$ is well-defined on
functions on $C_0\cap C_1$  for $\xi$ tangent to
$C_0\cap C_1$, and we set
\[
C^\infty_\epsilon(\underline{C_0\cap C_1})=\{f\in
C^\infty(C_0\cap C_1)[[\epsilon]] : A(E^\mu)f=0, \mu\in I_0^c
\cap I_1^c\}.
\]
As $A(E^\mu)f=E^\mu f+O(\epsilon)$ this condition 
reduces modulo $\epsilon$ to the condition that $f\in
C^\infty(\underline{C_0\cap C_1})$. 

If $\Gamma\in G_{k,2}$ we have two weights $w^0_\Gamma$
$w^1_\Gamma$ for the two module structures. The
weight $w^0_\Gamma$ is obtained by integrating over
the configuration space of $k$ distinct points in the first
quadrant, one point at the origin, associated to the
first vertex  $\bar 1$ and one point on the
positive real axis, associated to $\bar 2$, up to
dilations. The right $C^\infty_\epsilon(\underline {C_0})$\ndash module
structure is then defined by the product
\begin{equation}\label{e-Bradshaw}
\psi\star_0 f=\psi f+\sum_{k=1}^\infty
\frac{\epsilon^k}{k!}\sum_{\Gamma\in G_{k,2}}
w^0_\Gamma B_\Gamma(\psi,f),
\end{equation}
$\psi\in
C^\infty_\epsilon(\underline{C_0\cap C_1})$, 
$f\in C^\infty_\epsilon(\underline {C_0})$.
Similarly, the weights $w^1_\Gamma$ obtained by 
assigning $\bar 1$ to a point on the positive imaginary
axis and $\bar 0$ to the origin, we get the
left $C^\infty_\epsilon(\underline{ C_1})$\ndash module structure
\begin{equation}\label{e-PrimeMinister}
f\star_1\psi =\psi f+\sum_{k=1}^\infty
\frac{\epsilon^k}{k!}\sum_{\Gamma\in G_{k,2}}
w^1_\Gamma B_\Gamma(f,\psi),
\end{equation}
$\psi\in
C^\infty_\epsilon(\underline{C_0\cap C_1})$,
$f\in C^\infty_\epsilon(\underline {C_1})$.
Applying Stokes' theorem to this situation gives
our result:

\begin{Thm}\label{Luna}
Let the Poisson manifold 
$M$ be  an open subset of $\mathbb R^n$ containing
the origin and let $C_q$, $q=0,1$ be two coisotropic
submanifolds given by the equation $x^\mu=0, \mu\in I_q^c$.
Suppose that Assumption~\ref{a-Elisabeth} holds for
both $C_0$ and $C_1$. Then 
\begin{enumerate}
\item[(i)] The product $\star_0$ 
\eqref{e-Bradshaw} maps 
$C^\infty_\epsilon(\underline{C_0\cap C_1})
\otimes
C^\infty_\epsilon(\underline{C_0})
$ to $C^\infty_\epsilon(\underline{C_0\cap C_1})$ and
is a right $C^\infty_\epsilon(\underline{C_0})$\ndash module
structure.
\item[(ii)] The product $\star_1$ 
\eqref{e-PrimeMinister} maps 
$C^\infty_\epsilon(\underline{C_1})
\otimes C^\infty_\epsilon(\underline{C_0\cap C_1})
$ to $C^\infty_\epsilon(\underline{C_0\cap C_1})$ and
is a left $C^\infty_\epsilon(\underline{C_1})$\ndash module
structure.
\item[(iii)] We have $(f\star_1\psi)\star_0 g=
f\star_1(\psi\star_0 g)$, i.e., the two module structures
combine to a bimodule structure.
\item[(iv)] The reduction modulo $\epsilon$ is a 
homomorphism of bimodules.
\end{enumerate}
\end{Thm}
An important special case is the case where 
$C_0=M$. Then $\underline M=M$ and Assumption~\ref{a-Elisabeth} is
satisfied. Moreover, the algebra
$C^\infty_\epsilon(M)$
is the Kontsevich deformation of the product 
on $C^\infty(M)$ and 
$C^\infty_\epsilon(\underline{C_0\cap C_1})$
is $C^\infty(C_1)[[\epsilon]]$. In this way we get
a $C_\epsilon^\infty(M)^\mathrm{op}\otimes
C^\infty_\epsilon(\underline C)$\ndash module 
structure on $C^\infty(C)[[\epsilon]]$
for a coisotropic $C$ obeying Assumption~\ref{a-Elisabeth}.

\section{Examples}\label{exa}
Here we discuss some cases where Assumption~\ref{a-Elisabeth}
 is satisfied.
In all cases we assume that $M$ is an open subset of
$\mathbb R^n$ and that the coisotropic submanifolds are
coordinate subspaces, as in the previous Section.

\subsection{Codimension one} If $C\subset M$ is any coisotropic hyperplane,
Assumption~\ref{a-Elisabeth}
 is satisfied since the conormal bundle is
one-dimensional and $F$ is a skew-symmetric bilinear form.
So for each coisotropic hyperplane $C$ we obtain an algebra
$C^\infty_\epsilon(\underline C)$
quantizing the algebra of invariant functions on $C$, a
$C^\infty_\epsilon(M)^\mathrm{op}\otimes
C^\infty_\epsilon(\underline C)$\ndash module $C^\infty_\epsilon(C)$
and, for each pair of coisotropic hyperplanes $C_0$, $C_1$
a bimodule $C^\infty(\underline{C_0\cap C_1})$.

\subsection{Constant case} Let
 $M=\mathbb R^{n}$ with constant Poisson structure
and let 
$C$  be a coisotropic subspace. In this case Assumption~\ref{a-Elisabeth} 
is trivially satisfied as $F$ 
involves derivatives of $\pi^{ij}$. Also the condition
$A(E^\mu)f=0$ reduces to $E^\mu f=0$ so that 
$C^\infty_\epsilon(\underline C)
=C^\infty(\underline C)[[\epsilon]]$.
For example, consider the case of the standard symplectic
structure on $\mathbb{R}^{2m}$. Lagrangian subspaces are
coisotropic, with characteristic foliation consisting of one
leaf.
Thus $C^\infty_\epsilon(\underline C)$ is the one-dimensional $\mathbb R[[\epsilon]]$ free module of constant functions. Taking $C_0=M$ and $C_1=C$ we get a trivial left
module structure and $C_\epsilon^\infty(C)=C^\infty(C)[[\epsilon]]$ is a right $C^\infty_\epsilon(M)$\ndash module.
It is a formal version of the space of states in quantum
mechanics.
\subsection{Linear case} Let $\mathfrak g$ be a Lie algebra and
$M=\mathfrak g^*$ with Kostant--Kirillov Poisson structure.
The annihilator
$\mathfrak h^\perp$ of some Lie subalgebras $\mathfrak h$
is then
a coisotropic subspace of $M$ (see Example~\fullref{phoenix}).
It can be shown that Assumption~\ref{a-Elisabeth} is satisfied in
this case. 
As the Poisson structure is linear we may replace
smooth functions by polynomial functions. Our construction
gives then a quantization of 
$S(\mathfrak g/\mathfrak h)$ as an
 $S(\mathfrak g)$--$S(\mathfrak g/\mathfrak h)^\mathfrak h$\ndash bimodule. The quantization of $S(\mathfrak g)$ is the Kontsevich
deformation quantization $U=S_\epsilon(\mathfrak g)$. It is 
isomorphic to the universal enveloping algebra of $\mathfrak g$
with bracket $\epsilon[\ ,\ ]$ over $\mathbb R[[\epsilon]]$.
The quantization of $S(\mathfrak g/\mathfrak h)^\mathfrak h$
is an algebra 
$S_\epsilon(\mathfrak g/\mathfrak h)^\mathfrak h$.
In general 
$S_\epsilon(\mathfrak g/\mathfrak h)^\mathfrak h
/\epsilon S_\epsilon(\mathfrak g/\mathfrak h)^\mathfrak h$ is not $S(\mathfrak g/\mathfrak h)^\mathfrak h$,
so we do not have a deformation quantization in general.
We do, however, in the reductive case: 
\begin{Thm}\label{t-MrsWhitbread}
Suppose $\mathfrak h$ is a Lie subalgebra of a finite dimensional Lie 
algebra $\mathfrak g$. Assume that $\mathfrak h$ admits an
$\mathrm{ad}(\mathfrak h)$\ndash invariant complement. 
Then the algebra $S_\epsilon(\mathfrak g/\mathfrak h)^\mathfrak h$
is isomorphic to $S(\mathfrak g/\mathfrak h)^\mathfrak h[[\epsilon]]$ as
an $\mathbb R[[\epsilon]]$\ndash module. The products
$\star_0$, $\star_1$ of subsection~\ref{s-44}
define a $U^\mathrm{op}\otimes S_\epsilon(\mathfrak g/\mathfrak h)^\mathfrak h $\ndash %
module structure on the space 
$S(\mathfrak g/\mathfrak h)[[\epsilon]]$ of functions on
$\mathfrak h^\perp$.
\end{Thm}

A very particular case where the assumptions of the Theorem are satisfied is $\frh=\frg$.
In this case, $\frh^\perp$ is the origin of $\frg$, a zero of the Kostant--Kirillov
Poisson structure. The construction yields then a $U^\mathrm{op}$\ndash module structure
on $\bbR[[\epsilon]]$, that is, a character of the quantum algebra that deforms evaluation
at zero.

\section{Formality with submanifolds}\label{form}
The integrals over configuration spaces and the (bi)differential operators considered above can be generalized to the more 
general setting of multivector fields and multidifferential operators.
In the absence of branes,
Stokes' theorem gives then identities which in \cite{Kon} 
are formulated as the existence
of an $L_\infty$\ndash quasiisomorphism from 
the differential graded Lie algebra (DGLA) 
of multivector fields on 
$\mathbb R^{n}$ and the DGLA of multidifferential
operators. This is the local part of 
Kontsevich's formality theorem, 
one of whose important applications is the globalization
of the star product \cite{Kon,CaFeTo}.

In the presence of submanifolds, one should expect a similar
theorem to hold.
If $C\subset M$ is a submanifold, 
it is natural to introduce the DGLA $\mathcal V(M,C)=
\oplus_{j\geq{-1}}V^j(M,C)$
of {\em relative multivector fields}. The space $V^j(M,C)$
 consists of 
multivector fields $\pi\in\Gamma(M,\land^{j+1} TM)$ 
whose restriction
to $C$ vanishes on $\land N^*C$. The Schouten--Nijenhuis bracket
restricts to a bracket on relative multivector
fields, which therefore form a DGLA (with trivial differential).
On the other hand, let $A(M,C)=\Gamma(C,\land NC)$
be the graded commutative algebra of sections of the
exterior algebra of the normal bundle $NC=T_CM/TC$.
The Hochschild complex $C(A,A)=\oplus_j \mathrm{Hom}_\mathbb R
(A^{\otimes j},A)$ of the  graded Lie algebra $A=A(M,C)$
is then a DGLA with respect to the Hochschild differential
and the Gerstenhaber bracket. We then define
the DGLA of {\em relative multidifferential operators} $\mathcal D(M,C)$
to be the subalgebra of $C(A,A)$ consisting of 
multidifferential operators.

In the case where $M$ is an open subset of 
$\mathbb R^n$ and $C\subset M$ is given by equations 
$x^\mu=0$, $\mu=m+1,\dots n$, the Feynman rules
described in Sect.~\ref{quant} give rise to
linear maps
\begin{equation}\label{e-uk}
U_{k}\colon\land^k\mathcal V(M,C)
\to \mathcal D(M,C)[1-k],
\end{equation}
defined as a sum over all admissible graphs with
$k$ vertices of the first type as in
Sect.~\ref{quant} but with arbitrary valences and
number of vertices of the second type from which
wavy edges are allowed to emerge.

\begin{figure}[b]
\begin{picture}(0,150) (150,0)\scalebox{.7}
{\includegraphics{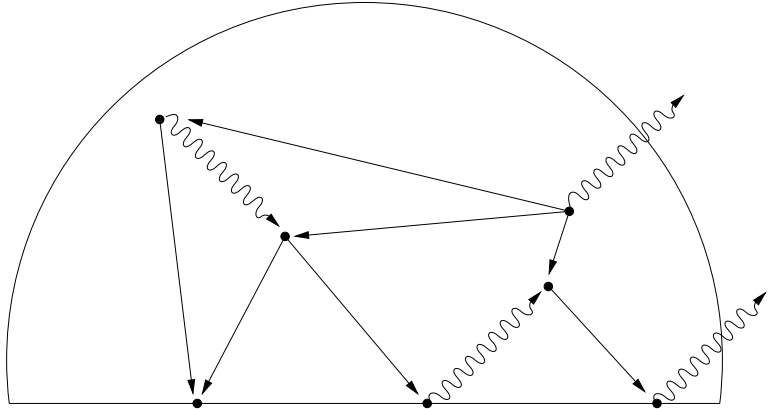}}
\end{picture}
\caption{A graph contributing to $U_4$}\label{fig2}
\end{figure}

\begin{Thm}\label{t-MissKilman} Let $C\subset M\subset
\mathbb R^n$ be the subset of the open set $M$ given
by equations $x^\mu=0,\mu=m+1,\dots,n$.
Then the maps $U_k$ are the Taylor coefficients of
an $L_\infty$\ndash morphism 
$U\colon\mathcal V(M,C)\to \mathcal D(M,C)$
\end{Thm}

In the case $M=C\subset \mathbb R^n$, 
$U$ reduces to Kontsevich's
$L_\infty$\ndash quasiisomorphism 
$\mathcal V(M)\to \mathcal D(M)$.

The DGLA $\mathcal D(M,C)$ is considered in disguise in \cite{OP} where it
is conjectured to be formal. 

 If we evaluate the maps $U_{k}$
on a Poisson
bivector field $\pi$, with $C$ coisotropic,
we recover the objects discussed
in subsection~\ref{s-43}: the solution 
 $U(\pi)=\sum \epsilon^k U_{k}(\pi,\dots,\pi)/k!$
of the Maurer--Cartan equation
in $D(M,C)[[\epsilon]]$ restricted to $C^\infty(C)\subset
A(M,C)$ has components of degree at most 2 in $NC$:
$U(\pi)=B(\pi)+\sum A(E^\mu)\de_\mu
+\frac12 \sum F(E^\mu,E^\nu)\de_\mu\land\de_\nu$. The star-product is $f\star g
=fg+B(\pi)f\otimes g$. 

An application of this generalized $L_\infty$\ndash morphism should be
the globalization of the deformation quantization of the reduced phase space
of a coisotropic submanifold. An extension of Thm.~\ref{t-MissKilman} to the case of two 
intersecting submanifolds should give a globalization
of the bimodule structure described in Thm.~\ref{Luna}
(though in general obstructions should be expected).

Let us add that we have considered here only the perturbative
part of the sigma model, namely the expansion of the 
path integral around a trivial classical solution.
The general case should lead to a generalization of the
Fukaya $A_\infty$\ndash category.

\subsection*{Note added} The $L_\infty$\ndash morphism of Thm.~\ref{t-MissKilman} can actually be defined on the Lie algebra $\mathcal V(M)$ of all multivector fields on $M$, not just on $\mathcal V(M,C)$. It induces \cite{CFcoisoclass} an $L_\infty$\ndash quasiisomorphism $\mathcal V(M)/I_C\to \mathcal D(M,C)$ on the quotient of $\mathcal V(M)$ by the Lie ideal $I_C$ consisting of multivector fields with vanishing Taylor expansion at each point of $C$. This quotient may be thought of
as the Lie algebra of multivector fields in a formal neighborhood of $C$.

\end{document}